    \documentclass[reqno]{conm-p-l}
    \usepackage{amsfonts}
    \usepackage{amssymb}

    \theoremstyle{plain}
   \newtheorem{thm}{Theorem}
   
   \newtheorem{lem}[thm]{Lemma}

   \theoremstyle{definition}

   \theoremstyle{remark}
   \newtheorem{rem}[thm]{{\it Remark}}

   \def\tlabel{\label}

\DeclareMathOperator{\okr}{{\stackrel{{\scriptscriptstyle{df}}}{=}}}
\DeclareMathOperator{\D}{d\!} \DeclareMathOperator{\SP}{supp}
\DeclareMathOperator{\CLOCONV}{cloconv}
 \DeclareMathOperator{\I}{i}
\DeclareMathOperator{\lin}{lin}

\DeclareMathOperator{\clo}{clo}
\def\ccc{\mathcal C}
\def\ddc{\mathcal D}

\def\hhc{\mathcal H}
\def\kkc{\mathcal K}
\def\mmc{\mathcal M}
\def\ccb{\mathbb C}
\def\nnb{\mathbb N}
\def\rrb{\mathbb R}
\def\Le{\leqslant}
\def\Ge{\geqslant}
\def\dz#1{\mathcal D({#1})}

\def\funkc#1#2#3#4#5{#1\colon#2\ni#3\mapsto#4\in#5}
\def\is#1#2{\langle#1,#2\rangle}

\def\liczp#1{{${#1}^{\text {\rm o}}$}}
\def\ulamek#1#2{\mbox{\normalfont$\frac{#1}{#2}$}}
\def\zb#1#2{\{{#1}\colon\ {#2}\}}

\begin{document}
   \title[Moments from their very truncations ]{Moments from their very truncations}
   \author[F.H. Szafraniec]{Franciszek Hugon Szafraniec}
   \address{Instytut        Matematyki,         Uniwersytet
   Jagiello\'nski, ul. Reymonta 4, PL-30059 Krak\'ow}
   \email{Franciszek.Szafraniec@im.uj.edu.pl}
   %
   \dedicatory{Andrzej Lasota (1932\,--\,2006) in memoriam}
   \begin{abstract}
   It is known that positive definiteness is not enough for the
   multidimensional moment problem to be solved. We would like
   throw in to the garden of existing in this matter so far results
   one more,  a result which takes into considerations the utmost
   possible truncations.
    \end{abstract}
   \maketitle As we have already pointed out positive definiteness
   is not sufficient for a multisequence to be a moment one,
   neither in the case of real moment problem in more than one
   variable, nor for a complex one or any complex dimension; for
   the previous one we recommend the cult paper of Fuglede
   \cite{fu}, for the second mentioned \cite{st-sz2} can be
   regarded as a source of information. Replacing it by
   solvability of a kind of truncations we gain necessary and
   sufficient conditions for the moment problem to be settled. It
   is worthy to say that truncations in the multivariable
   development have been considered from diverse points of view;
   let us have \cite{raul}, \cite{raul1}, \cite{put}, \cite{xu1}
   or \cite{xu2} as a choice of references.
   \subsection*{1}
   Let $\mmc(X)$ stand for the space of all regular complex Borel
   measures on a locally compact space $X$ and let $\mmc_a(X)$ be
   the collection of all positive measures in $\mmc(X)$ such that
   $\mu(X)=a$. Consider $\mmc(X)$ with the
   $\sigma(\mmc(X),\ccc_{\rm b}(X))$ topology\,\footnote{\;This
   topology is called `weak' one in Probability and Stochastic
   Processes though it may sound oddly for people in Functional
   Analysis.}, where $\ccc_{\rm b}(X)$ is the Banach space of
   continuous and bounded functions on $X$ with the `sup' norm;
   the topology is determined by the duality
   \begin{equation*}
   (\mu,f)\mapsto\int_Xf\D\mu,\quad
\mu\in\mmc(X),\;f\in\mmc_{\rm b}(X).
   \end{equation*}
    One of the pleasant features of the $\sigma(\mmc(X),\ccc_{\rm
   b}(X))$ topology is that $\mmc_a(X)$ is stable under the
   closure, another is that it coincides on $\mmc_a(X)$ with the
   $*$--weak topology.

   Let $\varXi$ be a linear space with a seminorm $p$. Set
   \begin{equation*}
\mmc_\varXi\okr \prod\nolimits_{\xi\in\varXi}\mmc_\xi,
   \end{equation*}
   where $\mmc_\xi=\mmc(X)$ for every $\xi\in\varXi$. Endow
$\mmc_\varXi$ with the Tychonoff topology based on that of
$\mmc(X)$. Having $\{\mu_\xi\}_\xi\in\mmc_\varXi$ define
   \begin{equation}\label{19.4}
\mu_{\xi,\eta}\okr\ulamek 1
4(\mu_{\xi+\eta}-\mu_{\xi-\eta}+\I\!\mu_{\xi+\I\!
\eta}-\I\!\mu_{\xi-\I \!\eta}),\quad \xi,\eta\in\varXi.
   \end{equation}

   The following selection result is in\,\footnote{\;Actually it
is stated and proved there for the $*$--weak topology in
$\mmc(X)$, however the proof brings over {\em verbatim} to the
$\sigma(\mmc(X),\ccc_{\rm b}(X))$ topology case. As an immediate
consequence we can replace inequality in condition (iii) of
\cite{sesqui} by equality and this results in equality in (b) here
(notice that due to (ii) the measures involved are positive). Also
we replace norm by seminorm which is still acceptable due to our
general version of Jordan--von Neumann Theorem therein. The method
is flexible enough to tolerate all these changes.} \cite{sesqui}.

   \begin{thm} \tlabel{t18.1}
Let $\mmc\subset\mmc_\varXi$ be a nonempty set and let $p$ be a
seminorm satisfying the parellelogram low
   \begin{equation*}
   p(\xi+\eta)^2+p(\xi-\eta)^2=2(p(\xi)^2+p(\eta)^2),\quad
\xi,\eta\in\varXi.
   \end{equation*}
 Suppose $\{\mu_f\}_f\in\mmc$ implies
   \begin{enumerate}
   \item[(i)] $\{\mu_{\I\!\xi}\}_\xi\in\varXi$ as well as
$\{t^{-2}\mu_{t\xi}\}_\xi\in\mmc$ for $t\in\rrb\setminus\{0\}$ and
$\{\ulamek 12\mu_{\xi+\eta}+\ulamek
12\mu_{\xi-\eta}-\mu_\eta\}_\xi\in\mmc$ for $\eta\in\varXi$;
   \item[(ii)] $\mu_0=0$, $\mu_{\xi+\eta}+\mu_{\xi-\eta}-2\mu_\eta\Ge
 0$;
   \item[(iii)] $\mu_\xi(X)=p(\xi)^2$, $\xi\in\varXi$.
   \end{enumerate}
   Then there is $\{\mu_\xi\}_\xi\in\CLOCONV(\mmc)$ such that
   \begin{equation}\label{cc}
   \mu_\xi\Ge0,\quad
   \mu_{\xi+\eta}+\mu_{\xi-\eta}=2(\mu_\xi+\mu_\eta),\quad
   \mu_{z\xi}=|z|^2\mu_\xi,\quad \xi,\eta\in\varXi,\;z\in\ccb.
   \end{equation}
   Consequently,
    for every Borel subset $\sigma$ of $X$ the mapping, cf.
 \eqref{19.4},
   \begin{equation}\label{ccc}
{(\xi,\eta)}\to{\mu_{\xi,\eta}(\sigma)}
   \end{equation}
is a positive Hermitian bilinear {\rm(=}positive
sesquilinear{\rm)} form on $\varXi$
    and
    \begin{equation}\label{cd}
\mu_{\xi,\xi}(X)=p(\xi)^2,\quad \xi\in\varXi.
   \end{equation}
   \end{thm}
   The proof relies on Markoff--Kakutani fixed point theorem.
   \begin{rem} \tlabel{td}
   If $\{\mu_\xi\}_\xi$ is as in the conclusion of Theorem
\ref{t18.1} then
   \begin{equation*}
   |\mu_{\xi,\eta}(\sigma)|\Le p(\xi)p(\eta),\quad
\xi,\eta\in\varXi\quad \text{$\sigma$ a Borel subset of $X$}.
   \end{equation*}
   Indeed, by the Schwarz inequality applied to the mapping of
\eqref{ccc} and then by \eqref{cd} we have (apparently
$\mu_\xi=\mu_{\xi,\xi}$ due to \eqref{19.4} and \eqref{cc})
   \begin{align*}
   |\mu_{\xi,\eta}(\sigma)|^2\Le\mu_\xi(\sigma)\mu_\eta(\sigma)
\Le\mu_\xi(X)\mu_\eta(X)\Le p(\xi)^2p(\eta)^2.
   \end{align*}

   \end{rem}

   \subsection*{2}
   For
      $n=(n_1,\dots,n_d)\in\nnb^d$, $N=\{0,1,2,\dots\}$ here, and
      for $x=(x_1,\dots,x_d)\in\rrb$ or $\ccb$ we hold up the
      notation: $|n|\okr n_1+\cdots+n_d$ and $x^n\okr
      x_1^{n_1}\dots x_d^{n_d}$. Moreover, a bit perversely,
      $\infty\okr(\infty,\dots,\infty)$. Notation for the basic
      zero-one $d$--tuples is shortened to
   \begin{equation} \label{2}
   e_i\okr(\delta_{m,i})_{m=1}^d,\quad i=1,\dots,d.
   \end{equation}

   A $d$--sequence $(a_n)_{n=0}^\infty$,
$a_n=a_{n_1,\dots,n_d}$, is said to be ($d$--dimensional {\it
real}) {\it moment} one if there is a positive measure $\mu$
on $\rrb^d$ such that
   \begin{equation*}
a_n=\int\nolimits_{\rrb^d}x^n\mu(\D x),\quad n\in\nnb^d.
   \end{equation*}
   Likewise, $(c_{m,n})_{m,n=0}^\infty$ is said to be a
($d$--dimensional {\it complex}) {\it moment} $2d$--sequence
if there is a positive measure $\nu$ on $\ccb^d$ such that
   \begin{equation*}
c_{m,n}=\int\nolimits_{\ccb^d}z^m\bar z^n\nu(\D z),\quad m,n\in\nnb^d.
   \end{equation*}
   Apparently there are two definitions of {\it positive
definiteness:} for $(a_n)_{n=0}^\infty$ as well as for
$(c_{m,n})_{m,n=0}^\infty$. These multisequences have to satisfy
the following conditions
   \begin{equation}
\text{$\sum\nolimits_{m,n}a_{m+n}\xi_m\bar\xi_n\Ge0$ for any
finite sequence of $(\xi_n)_n\subset\ccb$} \tag{${\sf
PD}_\rrb$}
   \end{equation}
   or
   \begin{equation}
\text{$\sum\nolimits_{m,n,k,l}c_{m+l,n+k}
\xi_{m,n}\bar\xi_{k,l}\Ge0$ for any finite sequence of
$(\xi_{m,n})_{m,n}\subset\ccb$.} \tag{${\sf PD}_\ccb$}
   \end{equation}
   It is commonly known that positive definiteness is sufficient
   for a multisequence to be a moment one but it is
   \underline{not} necessary except the $1$--dimensional real
   case. Our goal here is to present necessary and sufficient
   conditions for the moment problem to be solved.

   Let $\varXi$ denotes the linear space of all
   $d$--sequences $(\xi_n)_{|n|\Ge0}$ of complex numbers
   which are zero but a finite number.

For $\xi=(\xi_n)_n\in\varXi$ define new multisequences
$(a_n^{\xi})_{n=0}^\infty$ and
$(c_{m,n}^{\xi})_{m,n=0}^\infty$ as follows
   \begin{equation}\label{21.1}
a_n^\xi\okr\sum_{k,l}a_{n+k+l}\xi_k\bar\xi_l,\quad
c_{m,n}^\xi\okr\sum_{k,l}c_{m+k,n+l}\xi_k\bar \xi_l,\quad
m,n\in\nnb.
   \end{equation}
   Let us set
   \begin{equation*}
   {\sf N}_1\okr\zb{k}{|k|\Le 1}.
   \end{equation*}

      \begin{thm} \tlabel{t19.1}
A $d$--sequence $(a_n)_{n=0}^\infty$ is a moment one if and only
if there is a family $\{\mu_\xi\}_{\xi\in\varXi}$ of positive
measures on $\rrb^d$ such that
   \begin{equation} \label{19.1}
a_n^\xi=\int\nolimits_{\rrb^d}x^n\mu_\xi(\D x),\quad n\in{\sf
N}_1\cup2{\sf N}_1
   \end{equation}
   and
   \begin{equation} \label{19.2}
\mu_0=0,\quad \mu_{\xi+\eta}+\mu_{\xi-\eta}-2\mu_\eta\Ge0,\quad
\xi,\eta\in \varXi.
   \end{equation}
   \end{thm}

   The twin result, concerning the complex moment problem is as follows
   \begin{thm} \tlabel{t19.9}
A $2d$--sequence $(c_{m,n})_{m,n=0}^\infty$ is a moment one if and
only if there is a family $\{\mu_\xi\}_{\xi\in\varXi}$ of positive
measures on $\ccb^d$ such that
   \begin{equation*}
   c_{m,n}^\xi=\int\nolimits_{\ccb^d}z^m\bar z^n\mu_\xi(\D
z),\quad m,n\in{\sf N}_1
   \end{equation*}
   and
   \begin{equation*}
\mu_0=0,\quad \mu_{\xi+\eta}+\mu_{\xi-\eta}-2\mu_\eta\Ge0,\quad
\xi,\eta\in \varXi.
   \end{equation*}
   \end{thm}

   There is a rather formal link between $2d$--dimensional
   real moment problem and $d$--dimensional complex one (see,
   Proposition 57 in \cite{st-sz2}) however on the operator
   level, which is pretty often used in proofs, it becomes
   fragile. Fortunately, the proof of Theorem \ref{t19.9}
   goes the same way as that of Theorem \ref{t19.1} otherwise
   a reader may elaborate the aforesaid link for him/herself.

      \subsection*{3}
      {\sc Proof of Theorem \ref{t19.1}}. Putting $n=0$ in
  \eqref{19.1} we get (${\sf PD}_\rrb$). This allows us to
  construct a reproducing kernel Hilbert space\,\footnote{\;Though
  elements of the RKHS approach can be traced on many occasions we
  would like to advertise here \cite{ksiazka}, at least for those
  who can read it.} $\hhc$, say, composed of sequences of complex
  numbers. More precisely, after defining $a_{(n)}\okr(a_{n+k})_k$
  the scalar product of $\hhc$ is given as
   \begin{equation*}
\is{a_{(m)}}{a_{(n)}}\okr a_{m+n},\quad m,n\in\nnb^d.
   \end{equation*}
   The space $\ddc\okr\lin\zb{a_{(n)}}{n\in\nnb^d}$ is dense
in $\hhc$ and
   \begin{equation} \label{1}
a_{m+n}^\xi=\is{\,\sum\nolimits_k\xi_ka_{(m+k)}}
{\sum\nolimits_l\xi_la_{(n+l)}},\quad m,n\in\nnb^d.
   \end{equation}
Define on $\varXi$ a seminorm $p$ by
   \begin{equation*}
   p(\xi)=\|\sum\nolimits_k\xi_ka_{(k)}\|,\quad
\xi=(\xi_k)_k\in\varXi
   \end{equation*}
   and by $\is {\,\cdot\,}-_p$ the related semi-inner product.
   \begin{rem} \tlabel{txx1}
   It is clear that $p(\xi'-\xi'')=0$ if and only if
$\xi'-\xi''\in\varDelta$ where
    \begin{equation*}
   \varDelta\okr\zb{(\xi_k)_k\in\varXi}{\sum_k\xi_ka_{m+k}=
   0\text{ for all }m\in\nnb^d}
   =\zb{(\xi_k)_k\in\varXi}{\sum_k\xi_ka_{(k)}= 0}.
   \end{equation*}
The set $\varXi$
    is apparently a linear subspace of $\varXi$ and so is
$\widehat\varXi\okr\varXi/\varDelta$. Moreover, $\widehat\varXi$
is a unitary space, the mapping
$$\funkc h{\widehat\varXi}{\hat\xi}{\sum\nolimits_k\xi_ka_{(k)}} \ddc$$ is
well defined and becomes a unitary operator between
$\widehat\varXi$ and $\ddc$

   \end{rem}

If $\mmc$ is defined by means of \underline{all} $\mu_\xi$,
   $\xi\in\varXi$, determined by \eqref{19.1} and \eqref{19.2},
   then there is a matter of direct verification to check that all
   the assumptions of Theorem \ref{t19.1} are fulfilled; in
   particular $\mmc$ is non empty (notice that \eqref{19.2} with
   $\xi=0$ forces the measures $\mu_\xi$ to be positive). Now it
   is a right time to make use of Theorem \ref{t18.1}. However,
   before doing this notice that $\mmc$ is already a convex set.
   So we have got a family $\{\mu_\xi\}_\xi$ of positive measure
   such that
   \begin{equation}\label{y1}
\{\mu_\xi\}_{\xi}\in\clo(\mmc).
   \end{equation}
    and \eqref{ccc} as well as \eqref{cd} hold true. In
particular, sesquilinearity in condition \eqref{ccc} of Theorem
\ref{t18.1} supported by Remark \ref{td} allows us to define the
family $\{\mu_f\}_{f\in\ddc}$ by $\mu_f\okr\mu_{h(\hat\xi)}$ which
is a well defined measure as long as $\xi\in
h^{-1}(f)=\hat\xi\in\widehat\varXi$, cf. Remark
\ref{txx1}\,\footnote{\;The reason we have had to make the
quotient operation at this stage is not because the proof would
not work. It is grounded upon incomparably deeper circumstances.
For the $1$--dimensional \underline{real} case people pretty often
neglect this assuming additionally that $\varDelta=\{0\}$. While
this does not cause too much pain (except aesthetical discomfort),
in the others it may generate great loss: the would-be moment
measures supported on real algebraic sets, regardless the moment
problem itself is real or complex, may be out of game -- those
important measures are just encoded in $\varDelta$, for much more
on this look at \cite{3czlony} and \cite{3czlonki}. The RKHS
approach shows sign of its might!}. This brings us back to the
Hilbert space $\hhc$ with \eqref{ccc}
 to be satisfied after replacing $\xi,\eta$ by $f,g$ and the
semi-norm in \eqref{cd} to be the norm of $\hhc$. Now by standard
means we extend $\mu_{f,g}$ to the whole of $\hhc$ and find a
semispectral measure $F$ in $\hhc$ such that
   \begin{equation*}
   \mu_{f,g}=\is{F(\,\cdot\,)f}g,\quad f,g\in\hhc.
   \end{equation*}

   With the shorthand notation \eqref{2} in mind define the
operators $A_i$ with domain $\dz{A^i}\okr\ddc$ by $A_ia_{(m)}\okr
a_{(m+e_i)}$. The operators $A_i$ are symmetric, $\ddc$ is
invariant for each of them and they commute pointwise on $\ddc$.
With $A^n\okr A_1^{n_1}\cdots A_d^{n_d}$ we have by \eqref{1}
   \begin{equation} \label{19.5}
a_n^{\xi}=\is{A^n\sum\nolimits_k\xi_ka_{(k)}}
{\sum\nolimits_l\xi_la_{(l)}}, \quad
n\in\nnb^d,\;\xi\in\varXi.
   \end{equation}
    Due to \eqref{y1}, after all those identifications, we can say
for any $f\in\ddc$ there is a net
$\{\mu_f^\alpha\}_\alpha\subset\mmc$ such that
   \begin{equation} \label{3}
   \int\nolimits_{\rrb^d}\varphi\D\mu_f^\alpha\to
\int\nolimits_{\rrb^d}\varphi\D\is{ F(\,\cdot\,)f}f, \quad
\varphi\in\ccc_{\rm b}(\rrb^d).
   \end{equation}
    Let $E$ be a spectral measure which the Naimark dilation
of $F$ living presumable in a larger space $\kkc$ and let
$B_i$ be defined as
   \begin{equation*}
\is{B_ix}{y}_\kkc\okr\int\nolimits_{\rrb^d}t^i\is{E(\D t)x}
y_\kkc, \quad x\in\dz B,\;y\in\hhc,\quad |i|=1
   \end{equation*}
   each with its maximal domain. We want to know $B_i$'s are
selfadjoint extensions of $A_i$'s. For this we use Lemma
\ref{t21.2} twice.

   First we show that $\ddc=\dz {A_i}\subset\dz {B_i}$,
$i=1,\dots,d$. For this the working part of condition \liczp 1 of
Lemma \ref{t21.2} applied to $\varPhi(t)\equiv|t_i|^2$ reads as
   \begin{align*}
   \int\nolimits_{\rrb^d}|t_i|^2\mu_f^\alpha (\D t)
=a^f_{2e_i}\quad \implies\quad
\int\nolimits_{\rrb^d}|t_i|^2\is{E(\D t)f}f_\kkc \Le a^f_{2e_i}
   \end{align*}
   which means that $f$ is in the domain of $B_i$. Now according
to the definition of $\mmc$ we have, by \liczp 2 of Lemma
\ref{t21.2},
   \begin{align*}
   \is{A_if}f-\is{B_if}f=\int\nolimits_{\rrb^d}x_i\D\mu_f^\alpha -
\int\nolimits_{\rrb^d}x_i\is{E(\D x)f}f=0
   \end{align*}
   as well as
   \begin{align*}
\is{A_if}{A_if}-\is{B_if}{B_if}
=\int\nolimits_{\rrb^d}|t_i|^2\mu_f^\alpha(\D t) -
\int\nolimits_{\rrb^d}|t_i|^2\is{E(\D t)f}f=0
   \end{align*}
    All this gives us $A_i\subset B_i$ for $i=1,\dots,d$, cf.
\cite[\S 5]{nagy}.

    Now, because $\ddc$ is invariant for every $A^n$ it is so also
for $B^n$. Thus, due to \eqref{19.5}, using multiplicativity
properties of spectral integrals (see, Theorem 4, p. 135 in
\cite{bir}),
   \begin{multline*}
   a_n=\is{A^na_{(0)}}{a_{(0)}}=\is{B^na_{(0)}}{a_{(0)}}=
\int\nolimits_{\rrb^d}x\is {E(\D x)a_{(0)}}{a_{(0)}}
=\int\nolimits_{\rrb^d}x^n\mu(\D x), \\ n\in\nnb^d
   \end{multline*}
   with $\mu=\is{F(\,\cdot\,)a_{(0)}}{a_{(0)}}$. Thus
$(a_n)_n$ is a $d$--dimensional moment $d$--sequence
according to our wish.

   The `only if' part is a matter of straightforward verification.
     \qed

   \subsection*{4}
   Let us prove the main ingredient of the above proof because it may
   be interesting and useful for itself. For a topological space $X$
    denote by $\ccc_{\rm{c}}(X)$ the space of all continuous
    complex functions with compact support.
    \begin{lem}\label{t21.2}
    Let $X$ be a locally compact Polish space. If $\phi$ is a
continuous complex function on $X$, $\{\mu_\alpha\}_\alpha$ is a
net in $\mmc_1(X)$ and $\mu$ is a positive measure. Consider the
limit
   \begin{equation} \label{x1}
   \lim\nolimits_\alpha
\int_X\varphi\D\mu_\alpha=\int_X\varphi\D\mu.
   \end{equation}
    \begin{enumerate}
   \item[\liczp 1]  If \eqref{x1} holds for every $\varphi\in\ccc_{\rm c}(X)$ and $\int_{X}
\phi\D\mu_\alpha\Le c$ uniformly in $\alpha$ then $\int_{X}
\phi\D\mu\Le c$ provided $\phi\Ge0$;
    \item[\liczp 2] if  \eqref{x1} holds for every $\varphi\in\ccc_{\rm
b}(X)$ and $\int_{X} \phi\D\mu_\alpha\rightarrow a$ then $\int_{X}
\phi\D\mu= a$ provided $\phi$ is such that $\int_{X}|\phi|^2
\D\mu_\alpha\Le c$ uniformly in $\alpha$.
    \end{enumerate}
    \end{lem}

   \begin{proof}
   Consider a sequence $(\varphi_k)_k\in\ccc_{\rm c}(X)$ such that
$0\Le\varphi_k\nearrow 1$ pointwise. Then the sequence $\{\SP
\varphi_k\}_k$ of compact sets nests $X$.

   Taking the limit passage of the left hand side of
   \begin{equation*}
   \int_X\varphi_k\phi\D\mu_\alpha \Le\int_X\phi\D\mu_\alpha\Le c
   \end{equation*}
first in $\alpha$ then in $k$ we come up to \liczp 1.

For \liczp 2 take $\varphi_k(z)=1$ and write
    \begin{multline*}
    |a-\int_{X}\phi \D\mu|\Le
|a-\int_{X}\phi\D\mu_\alpha|+|\int_{X}(1-\varphi_k) \phi\D
\mu_\alpha|\\+|\int_{X}\varphi_k\phi\D\,(\mu_\alpha-\mu)|+
|\int_{X}(\varphi_k -1)\phi\D\mu|.
    \end{multline*}
    The second term plays the most sensitive role so let us treat it as follows.

  By the theorem of Prokhorov  (\cite{par}, theorem 6.7, p.47 or
\cite{str}, p. 121)
for a given $\varepsilon$ there is a compact subset $K$ of
$X$ such that $\mu_\alpha(X\setminus K)<\varepsilon$ for
any $\alpha$. Now, pick up $k_0$ so that $\varphi_k=1$ on $K$ for $k>k_0$ and write
   \begin{align*}
    |\int_{X}(1-\varphi_k) \phi\D \mu_\alpha| &\Le
    |\int_K(1-\varphi_k)\phi\D\mu_\alpha|+|\int_ {X\setminus
K}(1-\varphi_k)\phi\D\mu_\alpha|\\
    &\Le |\int_ { K}(1-\varphi_k)\phi\D\mu_\alpha|\\& +
\sqrt{\int_ {X\setminus K}|1-\varphi_k|^2\D\mu_\alpha\int_
{X\setminus K}|\phi|^2\D\mu_\alpha}
    \\& =\sqrt{\varepsilon c}.
    \end{align*}
  Notice the evaluation holds for all $\alpha$'s uniformly in
$k>k_0$. Because of this we can start with evaluating the
forth term going beyond $k_0$, if necessary, and being backed
by the Schwarz inequality
   \begin{equation*}
   |\int_{X}(\varphi_k
-1)\phi\D\mu|^2\Le\int_X|\varphi_k-1|\D\mu\int_X|\phi|^2\D\mu
   \end{equation*}
   and \liczp 1, then fixing $k$ in the third make this, fix
$\alpha$ in the first and finally take the advantage of the
evaluation for the second.
   \end{proof}

   \subsection*{5} We admit that the machinery we have used
is pretty heavy. This is so because we have patterned our proof on
the content of \cite{sesqui} and that concerns operators. We may
have a hope it can be done in more direct way and this may be a
kind of challenge.

   At first glance it looks like our result is of different nature
than of \cite{sto}. Instead of solving all the truncated moment
problems coming from a given multisequence, as required in
\cite{sto}, we confine ourselves to the (family of) the very
initial truncations. As both approaches provide us with necessary
and sufficient conditions they stimulate a question of comparing
their usefulness. If one agrees a truncated moment problem should
be solved in finitely supported measures to make things easier,
our truncations lead to algebraic conditions of order at most 2;
however there is a family of them to be solved, all of them
subject to the constrains of the type \eqref{19.2}. We count on
the invitation to be accepted.

\subsubsection*{Aknowlegments} The main result of this paper was
mentioned on several occasions. It was put forward for the first
time already during the combined 1996 Iowa events: {\it Workshop
on Recent Developments in Moments and Operators} and {\it AMS
Special Session on Moments and Operators;}
 the last presentation was at the Edwardsville 2006 conference,
since the latter author's glimmer of hope to give it a matured
written form has materialized in this account.

 The author would like to acknowledge his appreciation of Jan
Stochel's remarks on the previous version of this paper.

   \bibliographystyle{amsplain}

\begin{thebibliography}{99}
      \bibitem{bir}  M.\,S. Birman, M.\,Z. Solomjak,
{\it Spectral theory of self-adjoint operators in Hilbert space},
D. Reidel Publishing Company, Dordrecht/Boston/Lancaster/Tokyo,
{\bf 1987}.
      \bibitem{3czlony} D. Cicho\'n, J. Stochel and F.\,H. Szafraniec,
   Three term recurrence relation modulo an ideal and
orthogonality of polynomials of several variables, {\it J. Approx.
Theory},  {\bf 134}\,(2005), 11--64.
      \bibitem{raul} R.\,E.\, Curto, L. Fialkow, Solution of the
truncated complex moment problem for flat data,
{\it Memoirs Amer. Math. Soc.}, {\bf 568}\,(1996).
   \bibitem{raul1} \bysame, A duality proof of Tchakaloff's theorem,
{\it J. Math. Anal. Appl.}, {\bf 269}\,(2002), 519--532.
    \bibitem{fu} {B. Fuglede}, {The multidimensional
moment problem}, {\it Expo. Math.}, {\bf 1}\,(1983), {47--65}.
   \bibitem{par} {K.\,R. Parthasarathy},  {\it Probability  measures
on metric spaces}, {Academic Press}, {New York--London}, {\bf
1962}.
      \bibitem{put} M. Putinar, A dilation theory approach to
cubature formulas, {\it Expo. Math.}, {\bf 15}\,(1997), 183--192.
       \bibitem{st-sz2} J. Stochel, F. H. Szafraniec, The complex
moment problem and subnormality: a polar decomposition approach,
{\it J. Funct. Anal.} {\bf 159}\,(1998), 432--491.
   \bibitem{sto} J. Stochel, Solving the truncated moment problem
solves the full moment problem, {\it Glasgow J. Math.}, {\bf
43}\,(2001), 335--341.
   \bibitem{str} {D. Strook}, {\it Probability Theory, an
analytic     view},     {Cambridge     Univ.      Press},
{Cambridge, UK}, {\bf 1994}.
      \bibitem{sesqui} F.\,H. Szafraniec,   Sesquilinear selection
of elementary spectral measures and subnormality, in {\it
Elementary Operators and Applications}, Proceedings, Blaubeuren
bei Ulm (Deutschland), June 9--12, 1991, ed. M.Mathieu, pp.
243--248, World Scientific, Singapore, {\bf 1992}.
   \bibitem{zelazko} \bysame,    Subnormality and cyclicity, Banach
Center Publications, {\bf 67}\,(2005), 349--356.
   \bibitem{3czlonki} \bysame,   Favard's theorem modulo an ideal,
{\it Oper. Theory Adv. Appl.}, {\bf 157}\,(2005), 301--310.
    \bibitem{ksiazka}\bysame, {\it Reproducing kernel Hilbert spaces}
   [in Polish], Wydawnictwo Uniwersytetu Jagie\-llo\'nskiego,
Krak\'ow, {\bf 2004}.
   \bibitem{nagy} B. Sz.--Nagy,     {\it Extensions   of
linear transformations in Hilbert space which extend beyond this
space}, Appendix to F.Riesz, B.Sz.-Nagy, {\it Functional
Analysis}, Ungar, New York, {\bf1960}
   \bibitem{xu1} Y. Xu,     Cubature formulae and polynomial ideals,
{\it Adv. Appl. Math.}, {\bf 23}\,(1999), 211--233.
   \bibitem{xu2}    \bysame, Constructing cubature formulae by the
method of reproducing kernel, {\it Numer. Math.}, {\bf 85}\,(2000),
155--173.
   \end{thebibliography}
   
   \end{document}